\documentclass[11pt]{article}    
\usepackage{amssymb,amsmath,amsthm,amscd}
\allowdisplaybreaks 

\begin{document}           
\newtheorem{proposition}{Proposition}[section]
\newtheorem{theorem}[proposition]{Theorem}
\newtheorem{definition}{Definition}[section]
\newtheorem{corollary}[proposition]{Corollary}
\newtheorem{Lemma}[proposition]{Lemma}
\newtheorem{conjecture}{Conjecture}

\newtheorem{Remark}{Remark}[section]
\newtheorem{example}{Example}[section]

\newcommand{\ds}{\displaystyle}
\renewcommand{\ni}{\noindent}
\newcommand{\pa}{\partial}
\newcommand{\Om}{\Omega}
\newcommand{\om}{\omega}
\newcommand{\va}{\varepsilon}
\newcommand{\var}{\varphi_{y,\la}}

\newcommand{\la}{\lambda}
\newcommand{\sik}{\sum_{i=1}^k}
\newcommand{\vov}{\Vert\omega\Vert}
\newcommand{\Umy}{U_{\mu_i,y^i}}
\newcommand{\lamns}{\lambda_n^{^{\scriptstyle\sigma}}}
\newcommand{\chiomn}{\chi_{_{\Omega_n}}}
\newcommand{\ullim}{\underline{\lim}}
\newcommand{\bsy}{\boldsymbol}
\newcommand{\mvb}{\mathversion{bold}}
\newcommand{\R}{{\mathbb R}}
\newcommand{\bR}{{\mathbb R}}
\newcommand{\bC}{\mathbb{C}}
\newcommand{\bE}{\mathbb{E}}
\newcommand{\bH}{\mathbb{H}}
\newcommand{\bP}{\mathbb{P}}
\newcommand{\cF}{\mathcal{F}}

\newcommand{\beq}{\begin{eqnarray*}}
\newcommand{\eeq}{\end{eqnarray*}}

\newcommand{\ben}{\begin{enumerate}}
\newcommand{\een}{\end{enumerate}}

\newcommand{\beqs}{\begin{eqnarray*}&\displaystyle}
\newcommand{\eeqs}{&\end{eqnarray*}}

\renewcommand{\theequation}{\thesection.\arabic{equation}}
\catcode`@=11 \@addtoreset{equation}{section} \catcode`@=12



\begin{center}
{\LARGE
   C\`adl\`ag curves of SLE   driven by L\'evy processes  \\
}
\end{center}

\begin{center} {  Qing-Yang Guan\footnote{ School of Mathematics,  Loughborough University,  UK.
                   Email: q.guan@lboro.ac.uk.
   The author is supported  partly   by EPSRC/D05379X/1 and NSFC 10501048. }   }
\end{center}

\begin{abstract} Schramm Loewner Evolutions (SLE)   are random increasing hulls defined
through the    Loewner equation driven by Brownian motion. It is
known that the increasing  hulls are  generated by continuous
curves. When the driving process
  is of the form $\sqrt{\kappa} B+\theta^{1/\alpha}S$ for
   a Brownian motion $B$ and a symmetric
  $\alpha$-stable  process $S$ with $\kappa$ not equal to $4$ and $8$,
   we prove that the corresponding  increasing hulls are generated by
 c\`adl\`ag curves. \\[12pt]
  \emph{AMS 2000 subject classifications: Primary 60D05,  Secondary 60H10,   60G52 \newline
  }
\end{abstract}

\section{Introduction}
Let $(B_t)_{t\geq0}$ be a standard Brownian motion and
$\mathbb{H}=\{x+iy\in\bC:y>
  0\}$. The
chordal SLE is a class of random increasing hulls $(K_t)_{t\ge 0}$
in the closed   upper half plane $\overline{\bH} $,  which are
defined through the following stochastic Loewner        equation
  \begin{align} \label{leq} \partial_tg_t(z)&=\frac{2}{g_t(z)-U_t},
  \qquad g_0(z)=z,\qquad z\in\overline{\bH},
  \ \ \ \
  \end{align}
with $U_t=\sqrt{\kappa}B_t$ for some constant $\kappa\geq 0$. For
$z\in \overline{\mathbb{H}}$, function  $g_{\cdot}(z)$  in
(\ref{leq}) can be solved  when $t\in  [0,\zeta(z))$, where $\zeta(z)=\sup
\{t\geq 0: g_s(z)-U_s\neq 0, \ 0\leq s\leq t\} $ for $z\in
\overline{\mathbb{H}}\setminus \{0\}$ and $\zeta(0)=0$. For each
$t\geq 0$, $K_t$ is defined by the set of `broken points'   before
time $t$, i.e., $K_t=\{z\in\overline{\bH}:\zeta(z)\le t\}$. In
Rohde, Schramm \cite{ROS} (Lawler, Schramm, Werner \cite{LSW} for
$\kappa=8$), it is proved that the  increasing hulls are generated
by continuous curves $(\gamma(t))_{t\geq 0}$ in $\overline{\bH}$.
For $\kappa=2$ and $\kappa=8$ these continuous curves  are
originally the motivation of Schramm \cite{S} for introducing  SLE
as the possible scaling limit of the two dimensional  loop-erased random walk and the uniform
spanning trees, respectively. By \cite{ROS} we
know  that the increasing hulls are identity to the continuous curve
for $0\leq \kappa\leq 4$, and in general  $ \mathbb{H}\setminus
K_t$ is the unbounded connected component of $\mathbb{H}\setminus
\gamma[0,t]$. Let $f_t=g^{-1}_t$ and $\hat{f}_t(z)=f_t(z+U_t)$.
Almost surely, the curve can be defined by
\begin{align}\label{curve}\gamma(t)=\lim_{y\downarrow0}\hat{f}_t(iy),\ \ \ \ \forall t\geq
0.
\end{align}
In \cite{ROS}, it is also proved  that $(\gamma_t) $ is a simple
curve for $0\leq \kappa\leq4$; a self-intersecting path for
$4<\kappa<8$; and a space-filling curve for $\kappa\geq 8$.

The aim of  this paper is to  study the  c\`adl\`ag  structure of
the stochastic Loewner evolution  driven by certain
  L\'{e}vy   processes.     Let $0<\alpha<2$ and
  $(S_t)_{t\geq 0}$ be a standard
symmetric $\alpha$-stable process, i.e.,  a L\'evy process such that
$\mathbb{E}[\exp{\{\lambda S_t\}}]=e^{-t|\lambda|^{\alpha}}$ (cf.
Bertoin \cite{B} and Sato \cite{SAT}). In Guan, Winkel \cite{GW},
for $\kappa,\theta\geq 0$, stochastic Loewner        equation  (\ref{leq}) is studied for driving
processes \begin{align}
\label{as}U_t=\sqrt{\kappa}B_t+\theta^{1/\alpha}S_t
\end{align}
on phase transition properties. See also an earlier paper Rushkin,
Oikonomou, Kadanoff and
 Gruzberg \cite{RO} on these properties,  and their recent paper \cite{RO1}  for
  some physical motivation
 of  this model such as branching polymers. For $(U_t)$ in (\ref{as}), we
 use
    right derivatives  in (\ref{leq})  and     define
    $g_t,  K_t,f_t,\hat{f}_t,\gamma_t$
similarly  as above with  $\zeta(z)=\sup \{t\geq 0: g_s(z)-U_s\neq
0,\  g_s(z)-U_{s-}\neq 0,\ 0\leq s\leq t\} $ for $z\in
\overline{\mathbb{H}}\setminus \{0\}$. We may also use these
notations for other driving processes or functions after special
remarks.

There are some theoretical  motivations  for the study    of the
stochastic Loewner equations driven by jump processes.  It is easy to see that the jumps
of a driven function  can shift the increasing points on the boundaries  of  the hulls
and create some  tree-like clusters.
Hence the jump processes  may be used to construct  some interesting   domains
with  tree-like boundaries through  Loewner's equation which can not be constructed  by a continuous  driven function. In fact, some interesting domains in
   complex analysis have   a tree-like boundary.
In \cite{BS}, Beliaev and Smirnov     mentioned that SLE driven by
L\'evy processes  may be candidates  of fractal domains with high
multi-fractal spectrum.

  When considering a `nice'  c\`adl\`ag  driven function (c\`adl\`ag: right
continuous with  left limit), we may expect that  the increasing
hull can be generated by a  c\`adl\`ag  curve.  However, compared
with   a continuous curve, c\`adl\`ag curves can be far more
complicated than it seems. To reflect    their    difference,  we
prove  in the last section that a comb space, which is not  a locally connected space,
can be generated by a   c\`adl\`ag  curve.  In complex analysis,
locally connected is an important property for the boundary of a
domain because it is equivalent to  the continuous extension of a
conformal map taking this     domain as its image  (cf. Theorem 2.1
Pommerenke \cite{CH}). This  shows that   we can not prove the
continuous extension stated  in Theorem 1.1 below  as the continuous case in \cite{ROS}.

 Due to the
jumps of the (symmetric) $\alpha$-stable process, the corresponding
hulls have tree-like structure with infinite branching on the
boundary.
 In \cite{GW}, inspired  by the results of  SLE,  we conjecture that  the increasing
hulls $(K_t) $ are generated by c\`adl\`ag curves.  As the jumps for
a typical path of the $\alpha$-stable process are  dense on
$[0,\infty)$, this conjecture can not be deduced from  the Brownian
SLE curve directly. However,   a function can still be defined   by
(\ref{curve})  provided the limit exists. Then we can   study the
property of $(\gamma_t) $ and its relation to  the increasing hulls.
For any function $f$ on $[0,t]$,  set $f[0,t]=\overline{\{f(s):
0\leq s\leq t\}}$.
 A  main result of this paper is the
following. 
\begin{theorem}\label{aaa}
Let  $0<\alpha<2$,  $\theta\geq 0$ and
$\kappa\in(0,4)\cup(4,8)\cup(8,\infty)$. Assume that  the driving process $(U_t)
 $ in    (\ref{leq}) is   defined by (\ref{as}).
  Then, almost surely, the limits in
$(\ref{curve})$ are well defined  and $(\gamma(t)) $   is a
c\`adl\`ag function. Moreover, almost surely, the conformal maps
$(f_t)_{t\geq0}$ extend to $\overline{\mathbb{H}}$ continuously and
$(K_t)  $ is generated by $(\gamma(t)) $, i.e.,
$H_t:=\mathbb{H}\setminus K_t$ is the unbounded connected component
of $\mathbb{H}\setminus \gamma[0,t]$ for   $t\geq0$.
\end{theorem}
 
The case for $\kappa=0$ in Theorem 1.1  is also given in   Chen,  Steffen \cite{CR}  recently by a different method. 
Additionally, they   
prove that the Hausdorff dimension of $K_t$ is equal to one.  For   general $\kappa$, the Hausdorff dimension  might 
not depend on $\theta $ which we have  not studied  in this paper.    

  The rest part
of this section are to introduce the methods in the proof of Theorem
\ref{aaa}  and  to prepare some other  notations which will be used
throughout the paper.  We first consider  the
c\'adl\'ag property  in Theorem 1.1 for  the `truncated' driving process
$(\sqrt{\kappa}B_t+\theta^{1/\alpha}S_{\delta,t})$, where
$(S_{\delta,t})$ is   the truncated $\alpha$-stable process for some $\delta>0 $:
$$S_{\delta,t}:=S_t -\sum_{0<s\leq t}\triangle S_sI_{|\triangle
S_s|\geq \delta},\ \ \ \ \Delta S_t=S_{t}-S_{t-},\
S_{t-}=\lim_{s\uparrow t}S_{s}.  $$ Before introducing the
truncated case in details,
  we      see how to use it to prove  Theorem
  1.1.

For $\delta>0$, define by $T_1$   the first time $|\triangle
S_t|>\delta$
  and  define by induction $T_{n+1}=\inf\{t> T_n: |\triangle S_t|>\delta\}$ for $n\geq 1$.
 By the  c\`adl\`ag   property of  $(S_t) $,  almost surely we have
 $\lim_{n\rightarrow \infty}T_n=\infty$ and $T_{n+1}>T_n>0$ for $ n\geq1$.
    Before
   stopping time $T_1$, as    processes $(S_{\delta,t})  $
   and
$(S_t) $ are the same, the c\`adl\`ag property    in Theorem 1.1
follows by  the truncated   case directly. For the  time between $T_n$
and $T_{n+1}$, we need a fact  that $f_{T_n}(\cdot)$ can   extend   to
$\overline{\mathbb{H}}$ continuously (see Lemma \ref{I61Idfff}). With
the help of this fact  one can check by (\ref{leq}) that
  $$\gamma_t=f_{T_n}(\beta_{t-T_n}),\ \ \ t\in[T_n,T_{n+1}),$$
  where
  $(\beta_t) $ is the c\'adl\'ag  curve defined by the following truncated driving process
$$ U_{T_n}+\sqrt{\kappa}B_{t+T_n}+\theta^{1/\alpha}S_{\delta,t+T_n}
-\sqrt{\kappa}B_{ T_n}-\theta^{1/\alpha}S_{\delta, T_n},\ \ \ t\geq
0.$$ Thus, due to the continuity of
 $f_{T_n}$ on $\mathbb{\overline{H}}$,  the   c\`adl\`ag property of $
{\gamma}_t $ on  $[T_n,T_{n+1})$ follows by that
 property of
 $  {\beta}_t $.
 For the continuous
extension of $f_{T_n}(\cdot)$, the proof   is similar to  the continuous
extension of $\hat{f}_t(iy)$  in \cite{ROS}. Briefly, the $x$
coordinate in $f_{T_n}(x+iy)$   takes the role of the $t$ coordinate
in $\hat{f}_t(iy)$. See Lemma
 \ref{lemma} and \ref{I61Idfff} for   details.  We can further prove
 that
the  continuous extension
 of $f_t$ holds for all $t\geq0$.
 The proof is based on  a right continuity  of $
\partial K_t$ followed by  the c\`adl\`ag property of $(\gamma_t)$ (see (\ref{conti})). In general this
right continuous property is not true even for a continuous driving
function (see a counterexample in Marshall and Rohde \cite{MO}).

For the proof of the truncated case, the main steps is similar to
the continuous case in Rohde and Schramm \cite{ROS} and  Lawler
\cite{Law}. However, due to the jumps we need to overcome some
different difficulties. In \cite{ROS}, to prove the convergence in
(\ref{curve}), a key step is to obtain some estimates for the
complex derivative of $\hat{f}_t(\cdot)$, denoted by
$\hat{f}'_t(\cdot)$. By Lemma 3.1 \cite{ROS}, this reduces to the
derivative estimates for the backward flow of (\ref{leq}).
  Let $(\widetilde{B}_t)_{t\geq0}$ and $(\widetilde{S}_{\delta,
t})_{t\geq 0}$ be another Brownian motion and truncated symmetric
$\alpha$-stable process, respectively. We assume that
$(B_t)_{t\geq0}$, $(S_t)_{t\geq 0}$, $(\widetilde{B}_t)_{t\geq0}$
and $(\widetilde{S}_{\delta, t})_{t\geq 0}$ are independent
  and define
\begin{align} \label{gfffgg}  V_{t}:=\sqrt{\kappa}
 {B}_{t}+\theta^{1/\alpha}  {S}_{\delta,t},\ \ \
V_{-t}:=\sqrt{\kappa} \widetilde{B}_{t}+\theta^{1/\alpha}
\widetilde{S}_{\delta,t},\ \ \ \forall \ t\geq 0.\end{align}  Next
we take $(V_t)_{t\in\R}$ as the driving process in (\ref{leq}) and
extend the solution of   (1.1)   to the  negative times. The
notations $g_t, K_t,f_t,\hat{f}_t$ will
  be   also used in this setting.  For $t< 0$, the derivative in
(\ref{leq})  is defined by the left derivative and  $g_t$ is a conformal
map from $\mathbb{H}$ onto a subset of $\mathbb{H}$. The following
lemma is a straightforward extension of Lemma 3.1 in \cite{ROS}
which can be proved by the  stochastic Loewner        equation.
\begin{Lemma}\label{E:harmonic}
For driving process $(V_t)$ and for all fixed $t\geq 0$, the map
$z\rightarrow g_{-t}(z)$ has the same distribution as the map
$z\rightarrow \hat{f}_t(z)-V_t$.
\end{Lemma}In \cite{ROS},   the  estimates of
$|g_{-t}'(z)|$ for  SLE are obtained   by moment
estimates  which in turn  are given by constructing martingales.
Next we introduce this method in the
context of this paper. Let $\hat{z}=\hat{x} +i\hat{y}  \in
\mathbb{H}$ and set for $t\geq 0$
\begin{align}\label{sured}
z(t)=(x(t), y(t)):=g_{-t}(\hat{z})-V_{-t},\ \ \
\psi(t)=\frac{\hat{y} }{y(t)}|g_{-t}'(\hat{z})|.
\end{align}
 Direct
calculations show, \begin{align}\label{1}
dx(t)=\frac{-2xdt}{x^2+y^2}-d V_{-t},\ \ \ dy(t)=\frac{
2ydt}{x^2+y^2},\ \ \ d\log \psi(t) =\frac{-4y^2dt}{(x^2+y^2)^2} .
\end{align}In \cite{ROS},    some
function $F$ and constant $\mu$ are constructed so     that  the process
\begin{align}
M_t=\psi(t)^\mu F(z(t)),\ \ \ t\geq 0
\end{align}
is a martingale with respect  to the natural filtration. By It\^{o}'s
formula and (\ref{1}), the drift of the process $(M_t) $ is
$$\int_0^t\psi(s)^\mu\Lambda F(z(s))ds,$$ where $\Lambda$ is the
following operator
\begin{align}\label{Plane:from cone::0f}
\Lambda  := \frac{-4\mu y^2}{(x^2+y^2)^2}
-\frac{2x}{x^2+y^2}\partial_{x} + \frac{2y}{x^2+y^2}\partial_{y}
+\frac{\kappa}{2}
\partial^2_{x} + \theta\Delta^{\alpha/2}_{x|\theta^{1/\alpha}\delta} .
\end{align}
Here $\Delta^{\alpha/2}_{x|\delta}$ is the generator of
$(S_{\delta,t})$ given in  $(\ref{090})$ below.   Thus     a martingale
 $(M_t) $ can be constructed by choosing  a harmonic function
$F$ of  $\Lambda$, i.e., $\Lambda F=0$. Let $b\in\R$. When $\theta=0
$, it is proved in \cite{ROS} that function
\begin{align}\label{fff} F(x,y):=(1+(x/y)^2)^by^\nu
 \end{align} is harmonic with respect
to $\Lambda$ by setting
\begin{align}\label{10}
   \displaystyle \mu:=2b+\kappa  b(1-b)/2,\ \ \ \nu:=4b +\kappa
b(1-2b)/2.\end{align}
 For general L\'evy processes, it may not be possible to find
explicit  nontrivial  harmonic functions. Instead, we look for
super-harmonic functions which is enough for our purpose. To this end we   consider the function
$F$   in (\ref{fff})  with   some
different    parameters (see (\ref{EQ:h_1:h_2})). To prove the
super harmonic property of $F$, we control
$\Delta^{\alpha/2}_{x|\delta}F$ by $\partial^2_x F$. This  is
achievable   when   $0\leq \kappa< 8$ with $b\in(1/2,1]$. When $
\kappa>8$ and $b\in (0,1/2)$,  this can be done only for  $|x|/y$
sufficiently large. For smaller $|x|/y$, we obtain a weaker estimate
of $\Delta^{\alpha/2}_{x|\delta}F$  which can be     controlled by
the first term in  (\ref{Plane:from cone::0f})
when choosing a bigger  $\mu$. See more details in
Section 2.    For $\kappa=4$,
the estimates in Proposition \ref{II} is   not strong enough to
prove the continuous extension property  (cf.  Remark 4.1).   It is  interesting to know 
whether      Theorem 1.1   holds or not for 
$\kappa=4,8$  as the   curves in these two cases  are in a critical situation.  

 Another crucial ingredient in proving the   SLE curves in
\cite{ROS} is the pathwise increment  estimates  of the  Brownian
motion. In \cite{Law}, the role of this property is distilled into a
determined version.   For a determined continuous driving  function
$(u_t) $, a condition in Lemma 4.32 \cite{Law} is
\begin{align}\label{hdk}
\sup_{k/2^{2n}<t\leq (k+1)/2^{2n}}|u_t-u_{k/2^n}|\leq
c\sqrt{n}/2^n,\ \ \ \ 0\leq k\leq 2^n-1,\ n\geq 1
\end{align}
for some constant $ c>0$.  Due to the jumps, one can not expect
a     property like (\ref{hdk})  for the $\alpha$-stable
process or  the truncated $\alpha$-stable process. In fact, it might
not be easy to prove  whether
 the following
  weak local increment property is true or not:
\begin{align}\label{1AB}
\sup_{0\leq t\leq  1}\limsup_{h\downarrow 0}
{|S_{t+h}(\omega)-S_t(\omega)|}/{\sqrt{h\ln h}}< c(\omega),\ \ \ \ \
\  \ \ a.s..
\end{align}
Here, we give    another form of regularity for the  $\alpha$-stable process which may have independent interest.
This includes  a
kind of  uniform increment estimate  on  each interval
$[k/2^{2n},(k+1)/2^{2n})$  for  a truncated
  stable process with jumps less than  $2^{-\beta n}$, where  $\beta>1$.  On the other hand,   for $\beta<2/\alpha$,  the number of the jumps bigger than
$2^{-\beta n}$ in each interval
$[k/2^{2n},(k+1)/2^{2n})$ can be controlled by   a   finite number. See Lemma \ref{yes}
and \ref{sure} for details. According  to these two properties,
Lemma \ref{lemma} gives some sufficient conditions on the existence
of a c\`adl\`ag curve for the Loewner evolution.
 
The structure of this paper is the following. We give the moments
estimates of    $\hat{f}_t'$ in Section 2 and prove the path
increments properties of the symmetric $\alpha$-stable process in
Section 3. In Section 4 we give   a  determined version for the
existence of a c\`adl\`ag   curve  and apply it to prove Theorem
1.1. The proof of  Theorem 1.2  and an example are given   in the last section.

\section{Derivative   estimates for the truncated case}
Let $0<\alpha<2$. The generator of the truncated $\alpha$-stable
  process $(S_{\delta,t})$ is  the truncated fractional
Laplacian
\begin{align}\label{090}
\Delta_{x|\delta}^{\alpha/2}f(x):=\lim_{\varepsilon\downarrow
0}\mathcal{A}(1,-\alpha)\int_{\{y:\varepsilon<|y-x|<\delta\}}\frac{f(y)-f(x)}{|x-y|^{1+\alpha}}dy,
\end{align}
provided the limit exists,  where  $ \mathcal{A}(1,-\alpha)=
 {\alpha2^{\alpha-1}{\pi^{ {-1}/{2}}\Gamma( {(1+\alpha)}/{2})}/\Gamma(1-
 {\alpha}/{2}})$ and $f$ is a real function on $\R$ (cf. \cite{GW}).
    The
following   lemma is to compare operator
$\Delta_{x|\delta}^{\alpha/2}$ with Laplacian on some functions
which will be used later.
\begin{Lemma}\label{gg}
 Let $0<\alpha<2$ and
$a>0$. Set for $b\in (0,1]$
$$f(x)=
(1+ax^2)^b .$$ Then there exists a positive constant
$ C_1= C_1(\alpha,b)$ such that for $c>0$
\begin{align}\label{EQ:h_1:h_2eee}
   |\Delta_{x|c}^{\alpha/2}f(x)|&< C_1c^{2-\alpha}f''(x), \ \ \ \ \ \  \ \ \ \ \   \ \ \ \ \   \
     \ b\in( 1/2,1], \ x\in \R;
   \\
 \label{EQ:h_1:hss_62ee6e1}\Delta_{x|c}^{\alpha/2}f(x)&< C_1c^{2-\alpha}|f''(x)|, \ \ \
  \ \ \ \   \ \ \ \ \     \ \ \ \ b\in (0,1/2),\
 |x|>2^2((1-2b)a)^{-1/2},\ \alpha\neq 2b;
 \\
 \label{EQ:h_1:hss_662eee}\Delta_{x|c}^{\alpha/2}f(x)&< C_1c^{2-\alpha }(1+|\log c|) |f''(x)|,
   \  \ b\in (0,1/2),\
 |x|>2^2((1-2b)a)^{-1/2},\ \alpha= 2b ;\\\label{new1}|\Delta_{x|c}^{\alpha/2}f(x)|&< C_1ac^{2-\alpha},
     \ \   \ \ \ \
    \ \ \ \     \ \ \ \ \  \ \ \ \ \   \ \  \   b\in (0,1],\ x\in \R.\end{align}

\end{Lemma}\noindent{\bf Proof}$\ $   Direct calculation shows that
\begin{align}\label{ggg}f''(x)=2ab(1+a(2b-1)x^2)(1+ax^2)^{b-2}.
\end{align}
First we prove $(\ref{EQ:h_1:h_2eee})$. Let  $b\in( 1/2,1]$. By
(\ref{ggg}),
\begin{align}\label{o}
2ab(2b-1) (1+ax^2)^{b-1} \leq f''(x)\leq 2ab (1+ax^2)^{b-1}.
\end{align} By (\ref{o}) we have
\begin{align}
&|\Delta_{x|c}^{\alpha/2}f(x)|\label{ooo}\nonumber\\
=&\mathcal{A}(1,-\alpha)\left|\int_{x-c}^{x+c}\frac{f(y)-f(x)-f'(x)(y-x)}{|x-y|^{1+\alpha}}dy\right |\nonumber\\
=&
\mathcal{A}(1,-\alpha)\left|\int_{x-c}^{x+c}\int_x^y\frac{f''(u)(y-u)}{|x-y|^{1+\alpha}}dudy \right |  \nonumber\\
 \leq&
\mathcal{A}(1,-\alpha)\int_{-c}^{c}\int_{-t}^t\frac{|f''(x+u)|}{|t|^{\alpha}}du
dt
 \\
\leq&
\mathcal{A}(1,-\alpha)(2b-1)^{-1}f''(x)\int_{-c}^{c}|t|^{-\alpha}\int_{-t}^t(\frac{1+ax^2}{1+a(x+u)^2})^{1-b}du
dt\nonumber\\ \leq&
2\mathcal{A}(1,-\alpha)(2b-1)^{-1}(2-\alpha)^{-1}f''(x)c^{2-\alpha}\sup_{t>0}\int_{-1}^{1}
(\frac{1+ax^2}{1+a(x+tu)^2})^{1-b}du .\nonumber \end{align} We claim that
\begin{align}\label{oooo}M:=\sup_{t>0,a>0,x\in \R}\int_{-1}^{1}
(\frac{1+ax^2}{1+a(x+tu)^2})^{1-b}du<\infty.\end{align}
When
$|x|\geq 2t$,  we have $M<\int_{-1}^{1}
(\frac{1+ax^2}{1+ax^2/4})^{1-b}du<2^{3-2b}$.
 When
$|x|<2t$ and $\sqrt{1/a}<2t$, by the assumption that  $b\in(1/2,1]$, we
have for $t>0$
\begin{align}\int_{-1}^{1}
(\frac{1+ax^2}{1+a(x+tu)^2})^{1-b}du=&\frac{1}{t}\int_{-t}^{t}
(\frac{1/a+x^2}{1/a+(x+u)^2})^{1-b}du\nonumber\\ <&\frac{1}{t}\int_{-t}^{t}
\frac{2(2t)^{2-2b}}{ |x+u|^{2-2b}}du\leq \frac{2^{3}}{2b-1}. \nonumber \end{align} Similarly,  we have
for
$|x|<2t$ and $\sqrt{1/a}\geq 2t$
\begin{align}\int_{-1}^{1}
(\frac{1+ax^2}{1+a(x+tu)^2})^{1-b}du\leq 2^{2-b},\ \ \ \ t>0. \nonumber \end{align}
Combing the facts above, we  get
(\ref{oooo}) and hence  (\ref{EQ:h_1:h_2eee}) is true.

Next we assume that $b\in (0,1/2)$ and $|x|>2^2((1-2b)a)^{-1/2}$.
Assume  also that $x>0$ by symmetry.
By (\ref{ggg}), we have for
$|y|\geq 2((1-2b)a)^{-1/2}$
\begin{align}\label{o1q}
ab(1-2b) (1+ay^2)^{b-1} \leq| f''(y)|\leq 2ab(1-2b) (1+ay^2)^{b-1}.
\end{align}
By   (\ref{o1q})  we have for $|y|\leq
x/2$\begin{align}\label{ii}|f''(x+y)|\leq 2^{2(2-b)}|f''(x)|.
\end{align}   Therefore, if in addition that $x>2c$, we have by
(\ref{ooo})
\begin{align}\label{oo1q}
|\Delta_{x|c}^{\alpha/2}f(x)|\leq \mathcal{A}(1,-\alpha)(2-\alpha
)^{-1} 2^{2(3-b)}c^{2-\alpha}|f''(x)|.\end{align}
 Noticing that for $y\in \{z:|z-x|\geq
x/2\}$ we have $|y|\leq 3|x-y|$  and $$1+9a|x-y|^2\leq (a|x-y|^2)
 \frac{4 }{ax^2}+9a|x-y|^2\leq  10 a|x-y|^2.$$Hence for $
 2^2((1-2b)a)^{-1/2}<x\leq 2c$,  we have by
(\ref{ii}) and $f>0$
\begin{align}\label{1111}
&\mathcal{A}(1,-\alpha)^{-1}\Delta_{x|c}^{\alpha/2}f(x)\nonumber\\\leq&
 \int_{x-c}^{x+c}\frac{ f(y)-f(x)}{|x-y|^{1+\alpha}}
I_{|y-x|\geq x/2}dy+
\int_{x-c}^{x+c}\int_x^y\frac{f''(u)(y-u)}{|x-y|^{1+\alpha}}du
I_{|y-x|<x/2}dy\nonumber\\\leq &\int_{x-c}^{x+c}\frac{
(1+9a|x-y|^2)^{b}}{|x-y|^{1+\alpha}} I_{|y-x|\geq x/2}dy+(2-\alpha
)^{-1}  2^{2(3-b)}c^{2-\alpha}|f''(x)| \nonumber\\
\leq &2(10a)^b\int_{0}^{x+c}\frac{1}{|x-y|^{1+\alpha-2b}}  I_{|y-x|\geq x/2}dy+(2-\alpha
)^{-1}  2^{2(3-b)}c^{2-\alpha}|f''(x)| . \end{align}
If  we further assume that $\alpha>2b$,   (\ref{o1q}) and  (\ref{1111}) give
\begin{align}
\label{1111w}
&\mathcal{A}(1,-\alpha)^{-1}\Delta_{x|c}^{\alpha/2}f(x)\nonumber\\
\leq& \frac{2\cdot
10^b}{\alpha-2b}2^{\alpha-2b}(ax^2)^{b}x^{-\alpha}+(2-\alpha )^{-1}
2^{2(3-b)}c^{2-\alpha}|f''(x)| \nonumber\\\leq&
c_1c^{2-\alpha}|f''(x)|
\end{align} for some constant $c_1=c_1(\alpha,b)$.
If $\alpha<2b$,   (\ref{o1q}) and (\ref{1111}) give
\begin{align}
\label{1111wd}
&\mathcal{A}(1,-\alpha)^{-1}\Delta_{x|c}^{\alpha/2}f(x)\nonumber\\
\leq& \frac{2\cdot
10^b}{2b-\alpha} a^bc^{2b-\alpha}+(2-\alpha )^{-1}
2^{2(3-b)}c^{2-\alpha}|f''(x)| \nonumber\\\leq&
c_2c^{2-\alpha}|f''(x)|
\end{align}
for some constant $c_2=c_2(\alpha,b)$.
Thus  $(\ref{EQ:h_1:hss_62ee6e1})$ is given by (\ref{oo1q}), (\ref{1111w}) and (\ref{1111wd}).
The proof of (\ref{EQ:h_1:hss_662eee}) is   similar to
$(\ref{EQ:h_1:hss_62ee6e1})$.
 By (\ref{ggg}), we see that
$|f''(x)|\leq 2ab$ for $b\in(0,1]$ and hence  $(\ref{new1})$ follows
by (\ref{ooo}).  \qed\medskip

Next we turn  to the  derivative estimates for the truncated case.
Recall that  for $t\in \R$,  $V_{t}$ is defined  by (\ref{gfffgg}).
In this section, we consider the stochastic Loewner equation
(\ref{leq}) with $(U_t)$ replaced by  $(V_t)_{t\in\R}$ and adopt  the
previous  notations $g_t, K_t,f_t,\hat{f}_t$ in this setting.
        Let
$\hat{z}=\hat{x} +i\hat{y}  \in \mathbb{H}$. Define $z(t)=(x(t),
y(t))$ and  $ \psi(t)$ for $t\geq 0$ by (\ref{sured}). For $u\geq
\log \hat{y}$, set
\begin{align}
\widetilde{x}(u)=x({T_u(\hat{z})}),\ \ \
\widetilde{y}(u)=y({T_u(\hat{z})}),\nonumber
\end{align}
where
\begin{align}
T_u=-\sup\{t\in \R:\ \ \mbox{\emph{Im}}(g_t(z))\geq e^u\ \}.
\end{align}
Notice that $T_u$ is less than infinity almost surely which can be proved by the recurrence 
and the Markov property  of t$(V_t)$. 
 Let $b $, $\kappa_1, \kappa_2, a'\geq 0$  and define
\begin{align}\label{EQ:h_1:h_2}
 \left\{ \displaystyle \begin{array}{ll}
   \displaystyle a:=2b+\kappa_1 b(1-b)/2,\ \ \ \ \ \ \ \ \ \lambda:=4b +\kappa_1
b(1-2b)/2 ,\ \ \ 0\leq \kappa<8,&\ \ b\in (1/2,1];
   \\
 a:=2b+\kappa_2 b(1-b)/2+a',\ \ \  \lambda:=4b +\kappa_2
b(1-2b)/2  ,\ \ \ \  \ \ \ \  \kappa> 8,&\ \ b\in (0,1/2].
\end{array}\right. \end{align}

In the following proposition  and  Section 3, we need some results
of  L\'evy processes, such as the construction of L\'evy processes
by Poisson random measures and the It\^o's formula of L\'evy
processes.  We refer to  Applebaum \cite{APP} and  \cite{B} for the
details.

\begin{proposition}\label{II}
Let $b\in (0,1/2)\cup(1/2,1]$, $\kappa_1, \kappa_2, a'\geq 0$  and
define $a,\lambda$ by $(\ref{EQ:h_1:h_2})$.
  Let $0<\hat{y}<1$ and set for $u>\log \hat{y}$
\begin{align}
F(\hat{z})=\hat{y}^a  e^{-au}e^{\lambda
u} \mathbb{E}\left[(1+\widetilde{x}(u)^2/\widetilde{y}(u)^2)^b
|g'_{T_u(\hat{z})}(\hat{z})|^a\right].
\end{align}
1)\ If  $0\leq \kappa<8$, $\theta\geq 0$, $b\in (1/2,1]$ and
$\kappa_1>\kappa$, then there exists a constant
$\delta_0=\delta_0(\theta,\kappa_1-\kappa,\alpha,b)\in (0,1)$ such
that for any $0< \delta<\delta_0$,
\begin{align}\label{O}
F(\hat{z})\leq (1+(\hat{x} /\hat{y} )^2)^b\hat{y}^\lambda,\ \ \
0<\delta<\delta_0,
\end{align}
where $\delta$ is the parameter in (\ref{gfffgg}).

\noindent 2) If  $\kappa>8$ and $b\in (0,1/2)$, then we have the
same inequality (\ref{O}) for some constant
$\delta_0=\delta_0(\theta,\kappa-\kappa_2,\alpha,b)\in (0,1)$
provided $\kappa-\kappa_2
>0$ small enough.

\end{proposition}\noindent{\bf Proof}$\ $ 1) We first prove the case \textbf{$0\leq \kappa<8$},
$ b\in (1/2,1]$. Let $F_1( {z})= (1+( x/ y)^2)^by^\lambda$ and set
for $t\geq 0$
$$ M_t=\psi(t)^aF_1(z(t)).
$$
Let $\widetilde{N}(dt,dx)$ be the     Poisson random
martingale measure on $[0,\infty)\times
(-\theta^{1/\alpha}\delta,\theta^{1/\alpha}\delta)$ for the
truncated stable process
$(\theta^{1/\alpha}\widetilde{S}_{\delta,t})$.   The Poisson
intensity function of  $\widetilde{N}(dt,dx)$ is
\begin{align}\label{asa}(\mathcal{A}(1,-\alpha)\theta/|x|^{1+\alpha})I_{|x|<\theta^{1/\alpha}\delta}I_{t\geq 0}.\end{align}By
(\ref{1}) and It\^{o}'s formula (cf.  \cite{APP}) we have for
$t>0$
\begin{align}\label{Plane:from cone::0}
M_{t}=&M_0-\sqrt{\kappa}\int^{t}_0\psi(s)^a \partial_{x} F_1(z(s))\
d\widetilde{B}_{s}\nonumber\\
-& \int_0^t\int_{-\theta^{1/\alpha}\delta}^{\theta^{1/\alpha}\delta}
\psi(s)^a\left(F_1(z(s-)+x)-F_1(z(s-))\right)\widetilde{N}(ds,dx)+\int_0^t
\psi(s)^a \Lambda F_1(z(s))\ ds,
\end{align}
where $\Lambda$ is defined by (\ref{Plane:from cone::0f}) with $\mu$
replaced by $a$. Direct calculation shows
\begin{align}\label{Plane:from cone::000}
H(z):=\Lambda F_1(z)= \frac{\kappa-\kappa_1}{2}
\partial^2_{x}F_1(z)+ \theta\Delta^{\alpha/2}_{x|\theta^{1/\alpha}\delta}F_1(z),\ \ \ y>0.
\end{align}
By  Lemma \ref{gg}, $\kappa<\kappa_1$ and the fact that
$\partial^2_{x}F_1(z)>0$, we can choose
$\delta_0=\delta_0(\theta,\kappa_1-\kappa,\alpha,b)\in (0,1\wedge
\theta^{-1/\alpha})$ small enough such that $H(z)\leq 0$ for
$0<\delta<\delta_0$. Thus by (\ref{Plane:from cone::0}), we have
\begin{align}\label{Plane:from cone::0ss11}
M_{t}\leq
&M_0-\sqrt{\kappa}\int^{t}_0M_{s}\frac{2bx(s)}{x(s)^2+y(s)^2}\
d\widetilde{B}_{s}\nonumber\\
-& \int_0^t\int_{-\theta^{1/\alpha}\delta}^{\theta^{1/\alpha}\delta}
M_s\left(F_1(z(s-)+x)-F_1(z(s-))\right)F_1(z(s))^{-1}\widetilde{N}(ds,dx).
\end{align}

 Define
the right hand side of (\ref{Plane:from cone::0ss11}) by
$\widetilde{M}_t$ which is a local martingale.  Let
$u=u(\hat{z},t):=\log [\mbox{\emph{Im}} g_{-t}(\hat{z})]$. By
(\ref{leq}), $\partial_t
u=2(x(t)^2+y(t)^2)^{-1}=2(\widetilde{x}(u)^2+\widetilde{y}(u)^2)^{-1}$.
For $M>0$, set $T=\inf\{t\geq 0: |M_t|>M\}$. By $b> 1/2$ and the
mean value principle, for   $|x|<1$ we have
$$
|((\widetilde{x}(s)+x)^2+\widetilde{y}(s)^2)^b-(\widetilde{x}(s)^2+\widetilde{y}(s)^2)^b|\leq
2b|x|(|\widetilde{x}(s)|+1)^2+\widetilde{y}(s)^2)^{b-1/2}.$$ Hence
by (\ref{asa}) and  It\^{o}'s isometry,   we have for
$0<\delta<\delta_0$
\begin{align}\label{Plane:from cone::01111}
&\mathbb{E}{\widetilde{M}}_{T_u\wedge T}^2\nonumber\\\leq
&3M_0^2+3\kappa\mathbb{E}\left[\int^{T_u\wedge
T}_0M_{s}^2\frac{(2bx(s))^2}{(x(s)^2+y(s)^2)^2}\
d {s}\right]\nonumber\\
+& 3\mathcal{A}(1,-\alpha)\theta\mathbb{E}\left[\int^{T_u\wedge
T}_0M_{s}^2\int_{-\theta^{1/\alpha}\delta}^{\theta^{1/\alpha}\delta}
\frac{\left(F_1(z(s)+x)-F_1(z(s))\right)^2F_1(z(s)
)^{-2}}{|x|^{1+\alpha}}\
dxds\right]\nonumber\\
\leq &3M_0^2+3\kappa\mathbb{E}\left[\int^{
u}_{\log\hat{y}}M_{T_s\wedge
T}^2\frac{2b^2\widetilde{x}(s)^2}{\widetilde{x}(s)^2+\widetilde{y}(s)^2}\
d
{s}\right]+\frac{3}{2}\mathcal{A}(1,-\alpha)\theta\cdot\nonumber\\
 \mathbb{E}&\left[\int^{ u}_{\log\hat{y}}M_{T_s\wedge T}^2\int_{-\theta^{1/\alpha}\delta}^{\theta^{1/\alpha}\delta}
\frac{\left(((\widetilde{x}(s)+x)^2+\widetilde{y}(s)^2)^b-(\widetilde{x}(s)^2+\widetilde{y}(s)^2)^b\right)^2
(\widetilde{x}(s)^2+\widetilde{y}(s)^2)^{1-2b}}{|x|^{1+\alpha}}\
dxds\right]\nonumber\\
\leq&3M_0^2+6b^2\kappa\mathbb{E}\left[\int^{
u}_{\log\hat{y}}M_{T_s\wedge T}^2 \ d {s}\right]+6
\mathcal{A}(1,-\alpha)b^2\theta\cdot\nonumber\\&
\mathbb{E}\left[\int^{ u}_{\log\hat{y}}M_{T_s\wedge
T}^2\int_{-\theta^{1/\alpha}\delta}^{\theta^{1/\alpha}\delta}
\frac{((|\widetilde{x}(s)|+1)^2+\widetilde{y}(s)^2)^{2b-1}(\widetilde{x}(s)^2+\widetilde{y}(s)^2)^{1-2b}}
{|x|^{\alpha-1}}\ dxds\right]. \end{align} Noticing that
$\theta^{1/\alpha}<1$  and considering $|\widetilde{x}(s)|<1$ and
$|\widetilde{x}(s)|\geq 1$ respectively,  the   term  on   line of (\ref{Plane:from cone::01111})
  is less
than
\begin{align} \label{q}& \frac{12b^2\theta\mathcal{A}(1,-\alpha)}{2-\alpha}
\left((4/\hat{y}^2+1)^{2b-1}+2^{4b-2}\right)\mathbb{E}\left[\int^{
u}_{\log\hat{y}}M_{T_s\wedge T}^2 ds\right]
\leq c \mathbb{E}\left[\int^{
u}_{\log\hat{y}}\widetilde{M}_{T_s\wedge T}^2 \ d {s}\right]
\end{align}
for some constant $c=c(\kappa,\theta,b,\alpha, \hat{y})$. By
(\ref{Plane:from cone::01111}) and (\ref{q}), we have for $u\geq
{\log\hat{y}}$
\begin{align}
&\mathbb{E}{\widetilde{M}}_{T_u\wedge T}^2\leq3
M_0^2+c\mathbb{E}\left[\int^{u}_{\log\hat{y}}{\widetilde{M}}_{T_s\wedge
T}^2 \ d {s}\right].
\end{align}
Hence by Gronwall's Lemma (cf. \cite{RW}),
$$\mathbb{E}\left[ \widetilde{M}_{T_u\wedge T}^2\right]\leq 3M_0^2e^{c(u-{\log\hat{y}})}, $$ which gives the finiteness of   $\mathbb{E}\left[
\widetilde{M}_{T_u }^2\right]$  by taking $M\rightarrow
\infty$. This proves that $ (\widetilde{M}_{T_u})_{u\geq \log
\hat{y}}$ is a martingale. By optional stopping theorem   and
(\ref{Plane:from cone::0ss11}) we complete the proof of this case.

2) Assume that \textbf{$  \kappa>8$}, $ b\in (0,1/2)$ and $a'>0$.
  We only prove the case
$\alpha\neq2b$ as the proof for the case $\alpha=2b$ is similar. By
(\ref{EQ:h_1:hss_62ee6e1}), $(\ref{new1})$ and (\ref{ggg}) we have
for $\theta^{1/\alpha}\delta<1$
\begin{align}\label{1yes1}
\Lambda F_1(z)=& -\frac{4a'y^2}{(x^2+y^2)^2}(1+( {x}/{y})^2)^b
y^{\lambda}+ \frac{\kappa-\kappa_2}{2}
\partial^2_{x}F_1(z)+
\theta\Delta^{\alpha/2}_{x|\theta^{1/\alpha}\delta}F_1(z)\nonumber\\
\leq& -a'(\frac{2y^2}{x^2+y^2})^2(1+( {x}/{y})^2)^b y^{\lambda-2}+
b(\kappa-\kappa_2) \left(1+(2b-1)(x/y)^2\right)(1+(
{x}/{y})^2)^{b-2} y^{\lambda-2}\nonumber\\
+& C_1\theta(\theta^{1/\alpha}\delta)^{2-\alpha}\left|1+(2b-1)(x/y)^2\right|(1+(
{x}/{y})^2)^{b-2} y^{\lambda-2}I_{|x|>
2^2(1-2b)^{-1/2}y}\nonumber\\
+& C_1\theta(\theta^{1/\alpha}\delta)^{2-\alpha}y^{\lambda-2}I_{|x|\leq
2^2(1-2b)^{-1/2}y}.
\end{align}
 We claim that there exist   $\delta_0\in(0,1\wedge \theta^{-1/\alpha})$ and
$\kappa-\kappa_2>0$ small enough such that   $\Lambda F_1(z)\leq 0$
for $\delta\in(0,\delta_0)$. By (\ref{1yes1}) and $a'>0$, we can
first choose $\delta$ and
  $\kappa-\kappa_2$   small enough so that $\Lambda F_1(z)$ is negative
  for $|x|\leq
2^2(1-2b)^{-1/2}y$. By $\kappa>\kappa_2$ and $ b\in (0,1/2)$, we can
further decrease $\delta$ so that $\Lambda F_1(z)$ is negative
 for $|x|>
2^2(1-2b)^{-1/2}y$.

  Next we   prove that $
(\widetilde{M}_{T_u})_{u\geq \log \hat{y}}$ defined above  is still
a martingale in this case. By $2b-1<0$, we have for
$|\widetilde{x}(s)|<2$
$$N_1:=\sup_{|x|<1}(|\widetilde{x}(s)+x|^2+\widetilde{y}(s)^2)^{2b-1}
(\widetilde{x}(s)^2+\widetilde{y}(s)^2)^{1-2b}\leq
(1+2^2\hat{y}^{-2})^{1-2b},$$   and for $|\widetilde{x}(s)|\geq2$
$$N_2=:\sup_{|x|<1}(|\widetilde{x}(s)+x|^2+\widetilde{y}(s)^2)^{2b-1}
(\widetilde{x}(s)^2+\widetilde{y}(s)^2)^{1-2b}\leq 2^{2-4b}.$$ By
  these two estimates,    (\ref{Plane:from cone::01111}) is
replaced by
\begin{align}
&\mathbb{E}{\widetilde{M}}_{T_u\wedge T}^2\nonumber\\ \leq&
2M_0^2+(2b)^2\kappa\mathbb{E}\left[\int^{
u}_{\log\hat{y}}M_{T_s\wedge T}^2 \ d {s}\right]+ \mathcal{A} (1,-\alpha)(2b)^2\theta\cdot\nonumber\\
  &\mathbb{E}\left[\int^{
u}_{\log\hat{y}}M_{T_s\wedge
T}^2\int_{-\theta^{1/\alpha}\delta}^{\theta^{1/\alpha}\delta}
\frac{\sup_{|x|<\theta^{1/\alpha}\delta}(|\widetilde{x}(s)+x|^2+\widetilde{y}(s)^2)^{2b-1}
(\widetilde{x}(s)^2+\widetilde{y}(s)^2)^{1-2b}} {|x|^{\alpha-1}}\
dxds\right]\nonumber\\ \leq&
2M_0^2+(2b)^2\kappa\mathbb{E}\left[\int^{
u}_{\log\hat{y}}M_{T_s\wedge T}^2 \ d
{s}\right]+2\mathcal{A}(1,-\alpha)(2b)^2\theta(N_1+N_2)(2-\alpha)^{-1}\mathbb{E}\left[\int^{
u}_{\log\hat{y}}M_{T_s\wedge T}^2\ ds\right], \nonumber\end{align}
which gives
 the martingale property of $
(\widetilde{M}_{T_u\wedge T})_{u\geq \log \hat{y}}$. Therefore we can prove
the
 conclusion by
following the   arguments    in case 1). \qed\medskip

\section{Path increments of the symmetric $\alpha$-stable process}

For any    c\`adl\`ag  function $(u_t)_{t\geq 0}$, denote
$u_{t-}=\lim_{s\uparrow t}u_{s}$ and $\Delta u_t=u_{t}-u_{t-}$.  For
$\delta>0$, write $S_{\delta,t}=S_t -\sum_{0<s\leq t}\triangle
S_sI_{|\triangle S_s|\geq \delta}$ and  call
$(S_{\delta,t})_{t\geq0}$ the truncated symmetric $\alpha$-stable
process  with jumps less than $\delta$.
\begin{Lemma}\label{yes}
Let $0<\alpha<2$ and $(S_t)_{t\geq 0}$ be the symmetric $\alpha
$-stable process. Then for $\beta<2/\alpha$ and
$L=[3(2-\alpha\beta)^{-1}]+1$,
\begin{align}\label{P2P}
\lim_{m\rightarrow \infty}P\{\bigcup_{n\geq m}\bigcup_{0\leq
j<2^{2n}-1}A_{n,j}\}=0,\end{align} where \begin{align}  A_{n,j}=\{
\mbox{there exist }\ \ (t_i)_{i=1}^L\in [j/2^{2n},(j+1)/2^{2n}) \
\mbox{such that}\ |\triangle S_{t_i} |>2^{-\beta n },\ for\
i=1,2,\cdots,L \}.\nonumber
\end{align}
\end{Lemma}\noindent{\bf Proof}$\ $ We know that the jumps of the
symmetric $\alpha$-stable process  is a  Poisson point process  with
intensity  $(\mathcal{A}(1,-\alpha) /{|x|^{1+\alpha}})I_{t\geq 0}$
on $ \R\times \R_+$ (cf. \cite{B}). Thus for $n$ big enough
\begin{align}P\{A_{n,j}\}=&\exp\{-\mathcal{A}(1,-\alpha)\alpha^{-1}2^{1-(2-\alpha\beta)n}\}\sum_{k=L}^
\infty(\mathcal{A}(1,-\alpha)\alpha^{-1}2^{1-(2-\alpha\beta)n})^k(k!)^{-1}\nonumber\\
\leq&
(\mathcal{A}(1,-\alpha)\alpha^{-1}2^{1-(2-\alpha\beta)n})^L\nonumber\\
\leq & (2\mathcal{A}(1,-\alpha)\alpha^{-1} )^L2^{-3n},
\end{align}
which leads to (\ref{P2P}).\qed

\begin{Lemma}\label{sure}
Let $(S_t)_{t\geq 0}$ be the symmetric $\alpha $-stable process and
let $\beta>1$.   Then
\begin{align}
\lim_{m\rightarrow \infty}P\{\bigcup_{n\geq m}\bigcup_{0\leq
j<2^{2n}-1}B_{n,j}\}=0,\end{align} where \begin{align}\label{P1P}
B_{n,j}=\{\sup_{j/2^{2n}<t\leq (j+1)/2^{2n}}|S_{2^{-\beta
n},t}-S_{2^{-\beta n},{j/2^{2n}}}|>2^{-n}\}.
\end{align}
\end{Lemma}\noindent{\bf Proof}$\ $  Let $c>0$  and let
$S_{c,t}$ be the truncated symmetric $\alpha$-stable process with
jumps less than $c$. For $k\geq 1$, we denote
$$
\mathcal{P}_k=\{(2k_1,2k_2,\cdots,2k_l):\ 1\leq k_1\leq k_2\leq
\cdots\leq k_l,\ \sum_{s=1}^lk_s=k\ \mbox{for some positive integer
}l\geq 2\}.
$$
For each $(2k_1,2k_2,\cdots,2k_l)\in \mathcal{P}_k$, we say that a
partition $(A_k)_{k=1}^l$ for $\{1,2,\cdots,2k\}$ is of type
$(2k_1,2k_2,\cdots,2k_l)$ if the cardinal numbers of $A_k,1\leq
k\leq l$, are $(2k_1,2k_2,\cdots,2k_l)$ without considering the
order. Denote by $C(2k_1,2k_2,\cdots,2k_l) $  the cardinal number of
the type $(2k_1,2k_2,\cdots,2k_l)$ partitions for  the set
$\{1,2,\cdots,2k\}$. Set $(\triangle_t)_{ t\geq 0}=(\triangle
S_{c,t})_{t\geq 0}$. Since   $(\triangle_t)_{ t\geq 0}$  is the
Poisson point process on $\R^+\times \R$ with intensity
   $(\mathcal{A}(1,-\alpha)/{|x|^{1+\alpha}})I_{0<|x|<c}\cdot$ $
I_{t\geq 0}$,    by expansion and taking limit we have for integer $k\geq 1$
\begin{align}\label{moment}
 &\mathbb{E}|S_{c,t}|^{2k}\nonumber\\=&\lim_{\varepsilon\downarrow 0}\mathbb{E}|\sum_{0<s\leq
 t}\triangle_sI_{|\triangle_s|>\varepsilon}|^{2k}\nonumber\\=&\lim_{\varepsilon\downarrow
 0}\lim_{m\rightarrow\infty}
 \mathbb{E}|\sum_{i=1}^m\sum_{0<s\leq
 t}\triangle_tI_{\varepsilon+(c-\varepsilon)(i-1)/m<|\triangle_t|\leq \varepsilon+(c-\varepsilon)i/m} |^{2k}\nonumber\\
 =&\mathcal{A}(1,-\alpha)\int_0^t\int_{-c}^c\frac{|x|^{2k}}{|x|^{1+\alpha}}dx
 ds\nonumber\\
 +&\sum_{(2k_1,2k_2,\cdots,2k_l)\in
 \mathcal{P}_k}C(2k_1,2k_2,\cdots,2k_l)\prod_{i=1}^l\mathcal{A}(1,-\alpha)\int_0^t\int_{-c}^c
 \frac{|x|^{2k_i}}{|x|^{1+\alpha}}dxds
 \nonumber\\
 =& \frac{{2\mathcal{A}(1,-\alpha) }}{{2k-\alpha }}c^{2k-\alpha}t
 +\sum_{(2k_1,2k_2,\cdots,2k_l)\in
 \mathcal{P}_k}C(2k_1,2k_2,\cdots,2k_l)
  t^l c^{2k
 -l\alpha}\prod_{i=1}^l
 \frac{{2\mathcal{A}(1,-\alpha) }}{{2k_i-\alpha }}.
\end{align}
 Therefore for $t=c^{2/\beta}$ and $0<c<1$, we have
\begin{align}\label{moment2}
 &\mathbb{E}|S_{c,t}|^{2k}\nonumber\\\leq&    \frac{{2\mathcal{A}(1,-\alpha) }}{{2k-\alpha }}c^{2k+2/\beta-\alpha }
 +\sum_{(2k_1,2k_2,\cdots,2k_l)\in
 \mathcal{P}_k}C(2k_1,2k_2,\cdots,2k_l)\frac{{2^l\mathcal{A}(1,-\alpha)^l }}{({2-\alpha })^l}
  c^{2k+l(2/\beta-\alpha)
  }\nonumber\\\leq &  c_1 \max_{1\leq l\leq k} c^{2k+l(2/\beta-\alpha) }
\end{align}
for some $c_1=c_1(\alpha,k) $. By taking    $c=2^{-\beta n}$ in
(\ref{moment2}), we obtain
\begin{align}\label{P1P}
P\{B_{n,j}\}=&P\{\sup_{0<t\leq 1/2^{2n}}| S_{1/2^{\beta
n},1/2^{2n}}|>2^{-n}\}\nonumber\\
\leq &2P\{|S_{1/2^{ \beta n},1/2^{2n}}|>2^{-n}\}\nonumber\\\leq&
2^{1+2kn}\mathbb{E}|S_{1/2^{\beta n},1/2^{2n}}|^{2k} \nonumber\\
\leq&c_1 2^{1+2kn} \max_{1\leq l\leq k} 2^{-(2\beta
k+l(2-\alpha\beta))n }.
\end{align}
Noticing that $\beta>1$, we can choose $k$ big enough so that $P\{B_{n,j}\}\leq 2^{-3n}$, which gives the conclusion. \qed

\section{    C\`adl\`ag  curves}
    Lemma 4.32 in \cite{Law} gives some sufficient
conditions on the existence of the
  continuous curve   for the Loewner evolution. By checking    these
  conditions holding   almost surely,
   the existence of the  Brownian SLE curve can be proved. Here we
   adopt this strategy.
 Motivated by the
  regularity of the  stable processes studied in the last section,
  we put some conditions on the  jumps and fluctuations of a   c\`adl\`ag
  driven function
  to   give an extension of    Lemma 4.32 \cite{Law}.
  Lattice decomposition of time-space is important   in
proving the existence of  the continuous
      curve. Here we  divide each rectangle
in a lattice  decomposition of $\{(t,y)\in [0,N]\times[0,1]\}$
according to  the jumps and fluctuations of a   driven function. The similar method
was also used in \cite{ROS} for the continuous case.

Let $(u_t)_{t\geq 0}$ be a real c\`adl\`ag function. Suppose that
$(g_t)_{t\geq 0}$ is the Loewner chain in (\ref{leq}) with
$U_t=u_t$. For convenience,   the previous notations
$g_t,f_t,\hat{f}_t$ and $ \gamma_t$ will be also used for this
determined setting.
  Set $V(t,y)=\hat{f}_t(iy)$ for
$t\geq 0$  and $y>0$. Let  $\beta>0$ and
denote  for each integer $j\geq 1$, $u_{j,t}=u_t-\sum_{0\leq s\leq t}\triangle u_sI_{|\triangle
u_s|>2^{-\beta j}}$. By differentiating  both sides of
$f_t(g_t(z))=z$ with respect to $t$, we have (cf. \cite{Law})
\begin{align}
\label{02}
\partial_tf_t(z)=&-\frac{2f_t'(z)}{z-u_t},
  \qquad  \\
  \label{20}\partial_tf_t'(z)=&
  -\frac{2f_t''(z)}{z-u_t}+\frac{2f_t'(z)}{(z-u_t)^2}.
\end{align}
For any conformal  map $f$ defined  on the unite disc  to itself
with $f(0)=0$ and $|f'(0)|=1$, we know that $|f''(0)|\leq 2$  by
Bieberbach. This implies that for any conformal map  $g$ from
$\mathbb{H}$ to $\mathbb{H}$,
\begin{align}\label{kj}
|g''(z)|\leq 2|g'(z)|/\mbox{\emph{Im}}(z),\ \ \ \ z\in \mathbb{H}.
\end{align}
The  following lemma is a  c\`adl\`ag version of Lemma 4.32 in
\cite{Law} and the proof is similar.
\begin{Lemma}\label{lemma1}
Let $\beta>1 $ and   $L, N$ be  positive integers. Let
$(r_j)_{j\geq 1}$ be a sequence of increasing positive numbers with
$ \sum_{j=1}^\infty 1/r_j<\infty $. Suppose there exist  constants
$c,j_0>0$ such that, for each integer $j\geq j_0$ and $0\leq k\leq
N2^{2j}-1$, there exist  two increasing sequences
$(t(j,k,l))_{l=1}^{L}$ and  $(s(j,k,l))_{l=1}^{r_j}$ both belonging to
$[k/2^{2j}, (k+1)/2^{2j})$ and satisfying
\begin{align}\label{new}
\{t: |\triangle u_t|&>2^{-\beta j},\ t\in [k/2^{2j},
(k+1)/2^{2j})\ \}\subseteq \{t(j,k,l): l=1,2,\cdots, L\}, \\
\label{nnn}
 &\ \sup_{0\leq
l\leq r_j}\sup_{s(j,k,l)\leq s,t< s(j,k,l+1)}|u_{j,t}-u_{j,s}|\leq
c{2^{-j}} ,
\end{align}
where $s({j,k,0})=k/2^{2j}$ and $ s({j,k,r_j+1})=(k+1)/2^{2j}$.
Suppose that for $j\geq j_0$ and $0\leq k\leq N2^{2j}-1$ we also
have
\begin{align}\label{987}
|\hat{f}_{t(j,k,l)}'(i2^{-j}-\triangle u_{t(j,k,l)}) |&\leq
2^{j}/r_j^2,\ \ \ \ \  \ l=1,\cdots,L;\\ \label{d111d}\ \
|\hat{f}_{s(j,k,m)}'(i2^{-j}) |&\leq 2^{j}/r_j^2 ,\ \ \
    \ \ m=0,1,\cdots,
  r_j+1.
 \end{align}\emph{}
\noindent Then $\gamma(t)$ in (\ref{leq}) is well defined and
 is a c\`adl\`ag function on $[0,N]$.
\end{Lemma}\noindent{\bf Proof}$\ $  Let $j\geq j_0, 0\leq k\leq N2^{2j}-1$, $1\leq l\leq L$
and define $t(j,k,L+1)=t(j,k+1,1)$ with convention that
$t(j,N2^{2j},1)=N$. We see that
 $$
[t(j,k,l),t(j,k,l+1))\subseteq [k/2^{2j}, N\wedge(k+2)/2^{2j}).$$
Hence there are at most $2(r_j+1)$ numbers from  $ \cup_{ 0\leq
k\leq N2^j-1}\{(s(j,k,l))_{l=0}^{r_j} \}$ belonging to   the
interval $ [t(j,k,l),t(j,k,l+1)) $. Denote them by $a_1\leq
a_2\cdots\leq a_J$ and define $a_0=t(j,k,l),a_{J+1}=t(j,k,l+1)$. By
(\ref{20}) and (\ref{kj}) we have $|\partial_t f_t'(z)|\leq  6
|f_t'(z)|/(\mbox{\emph{Im}}(z))^2$ and hence
\begin{align}\label{aba}
  \ \ | f_{t-s}'(z)|\leq \exp\{6s/(\mbox{\emph{Im}}(z))^2\} |f_t'(z)| ,\ \ 0<s<t.
\end{align}
By (\ref{987}) and (\ref{aba}), we have
\begin{align}\label{amma}
&| {f}_{t}'(2^{-j}i+ u_{t(j,k,l+1)-} )|\leq e^{12}2^{j }/r_j^2 ,\ \
\
 \mbox{for}\ \ t\in
[a_J,a_{J+1}).
\end{align}
  By Distortion Theorem (see e.g.   \cite{Law}\cite{CH}), there exists a constant $K$ such that for
any conformal map $f$ from  $ \mathbb{H} $ to $ \mathbb{H} $,
\begin{align}\label{011}|f'(w)|\leq  K^{(|z-w|/y)+1}|f'(z)|, \ \ \
\mbox{\mbox{\emph{Im}}}(z), {\mbox{ {\emph{Im}}}}(w)\geq y>0.
\end{align} Combing (\ref{nnn}),   (\ref{amma}) and (\ref{011}), we
have for
  $t\in [a_J,a_{J+1})$
  \begin{align}\label{a2b2a}
&| \hat{f}_t'(i2^{-j} )|=| {f}_{t }'(i2^{-j}+u_{t(j,k,l+1)-}
+(u_t-u_{t(j,k,l+1)-} )| \leq   K^{c+1}e^{12} 2^{j }/ r_j^2.
\end{align}
By (\ref{011}) and (\ref{a2b2a})     we get
\begin{align}\label{yy}
| \hat{f}_t'(iy)|\leq  K^{c+3}e^{12} 2^{j} /r_j^2,\ \ \ \
2^{-j-1}\leq y\leq 2^{-j},\ \ t\in [a_J,a_{J+1}) .
\end{align}
Similarly, by (\ref{d111d}) we can also prove
\begin{align}\label{y}
| \hat{f}_t'(iy)|\leq  K^{c+3}e^{12} 2^{j} /r_j^2,\ \ \ \
2^{-j-1}\leq y\leq 2^{-j},\ \ t\in [a_m,a_{m+1} ),\ \ 0\leq m\leq
J-1.
\end{align}
Noticing that the estimates in (\ref{y})  only depend  on $j$, we
have for any $t\in [0,N]$
\begin{align}\label{03} | \hat{f}_t(iy)-\hat{f}_t(i2^{-j})|\leq
 K^{c+3}e^{12} \sum_{i=j}^\infty r_i^{-2},\ \ \ 0<y<2^{-j},
\end{align}which converges to zero as
$j\rightarrow \infty$  by the condition $\sum_{j\geq1}1/r_j<\infty$.
This implies that $\gamma(t)$ in (\ref{curve}) is   well defined for
$t\in[0,N]$.

By (\ref{02}) we have $|\partial_tf_t(z)|\leq 2|f'_t(z)|y^{-1}$.
Applying  (\ref{011}) one can also check that
$|\hat{f}'_t(i2^{-j}+x)|\leq  K^{2(c+1)}e^{12} 2^{j} /r_j^2$ if
$|x|<c2^{-j}$ and $t\in[0,N]$. Therefore, by (\ref{nnn}) and
(\ref{a2b2a})-(\ref{03}), we have for $0<y,y_1<2^{-j}$ and $t,t+s\in
[a_m,a_{m+1})$, $0\leq m\leq J$,
 \begin{align}
& |\hat{f}_t(iy)-\hat{f}_{t+s}(iy_1)|\nonumber\\\leq&
|\hat{f}_t(iy)-\hat{f}_{t
}(i2^{-j})|+|\hat{f}_t(i2^{-j})-\hat{f}_{t+s}(i2^{-j})|+|\hat{f}_{t+s}(iy)-\hat{f}_{t+s}
(i2^{-j})|\nonumber\\\leq& 2\cdot K^{c+3}e^{12}\sum_{i=j}^\infty
r_i^{-2}
+|\hat{f}_t(i2^{-j})-\hat{f}_{t}(i2^{-j}+u_{t+s}-u_t)|\nonumber\\+&
|\hat{f}_t(i2^{-j}+u_{t+s}-u_t)-\hat{f}_{t+s}(i2^{-j} )|
\nonumber\\
\leq& 5\cdot K^{2c+3}e^{12}\sum_{i=j}^\infty   r_i^{-2}.\nonumber
 \end{align}
As $J\leq 2(r_j+1)$ and $r_i$ is increasing, this gives that for
$0<y,y_1<2^{-j}$ and $t_1,t_2\in [a_0,a_{J+1})$
 \begin{align} |
\hat{f}_t(iy)-\hat{f}_{t+s}(iy_1)| \leq &5\cdot
K^{2c+3}e^{12}\sum_{i=j}^\infty  2(r_j+1) r_i^{-2} \nonumber\\\leq&
10\cdot K^{2c+3}e^{12}\sum_{i=j}^\infty   r_i^{-1},\nonumber
 \end{align}
 which gives  the conclusion  by
 $\sum_{j\geq1}1/r_j<\infty$.\qed\medskip

  From the proof above we   see that $(\gamma_t)$
   is continuous on the continuous  point of $(u_t)$. Next we
   denote by $(\Omega,\mathcal{F},P)$ the probability space  that  all the   stochastic  processes
   are considered. In the following  proposition
 the driving process is
 $(V_t)_{t\in\R}$ defined by (\ref{gfffgg}) with $\delta\in
(0,\delta_0)$, where  $\delta_0$ is the positive constant specified
by  Proposition \ref{II}.
\begin{proposition}\label{oiu}
Let $\kappa\in [0,8)\cup(8,\infty)$. Let  $b\in (0,1/2)$ for
$\kappa>8$ and $b\in (1/2,1]$ for $\kappa<8$. Define     $a,\lambda$
  by (\ref{EQ:h_1:h_2}) and suppose that
$\kappa_1,\kappa_2,a',\delta,\delta_0$ satisfies all the assumptions
in Proposition \ref{II}. For $T,t_1,t_2, a_1,a_2>0$ with
$0<t_2<t_1<T$, define stopping times $(T_k)_{k\geq 0}$ and
$(S_k)_{k\geq 0}$ by induction
 \begin{align}
T_0=& t_1,\ \ \ \ \ \ \ \ \ \ \ \ T_k=\inf \{ t<T_{k-1}: |\triangle
V_{t } |>a_1 \}\vee t_2,\ \ \ \ \
\ \ \ \ \ \ \ k\geq 1, \nonumber\\
S_0=& t_1,\ \ \ \ \ \ \ \ \ \ \  \ S_k=\inf \{ t<S_{k-1}: |  B_{t
}-B_{S_{k-1}} |>a_2 \}\vee t_2,\ \ \ k\geq 1.\nonumber
 \end{align}   Then there
exist  constants $ C_2(\kappa,,\kappa_1,\kappa_2,\alpha,
\theta,b), C_3(\kappa,,\kappa_1,\kappa_2,\alpha, \theta,b)$ such that
the following estimates hold  for all  $y,\rho\in(0,1]$, $x\in \R$
and $k=0,1,2,\ldots$,
\begin{align}\label{cd}
P\{|\hat{f}_{T_k}'(x-\triangle V_{T_k}+iy)|>\rho y^{-1}\}\leq&  C_2
 (1+(x/y)^2)^{b}(y/\rho)^{\lambda}\vartheta(\rho,a-\lambda),\ \
 \  \\ \label{cde}
 P\{|\hat{f}_{S_k}'(x +iy)|>\rho
y^{-1}\}\leq &  C_3
 (1+(x/y)^2)^{b}(y/\rho)^{\lambda}\vartheta(\rho,a-\lambda),\ \ \
\end{align}
 where
\begin{align}
\vartheta(\rho,s)=\left\{ \displaystyle \begin{array}{ll}
   \displaystyle \rho^{-s},\  \ \ \ &s>0,
   \\
 \ds 1+|\log \rho|,\ \ \ \ &s=0,\\
 1,\ \ \ \ &s<0.
\end{array}\right.\nonumber \end{align}

\end{proposition}\noindent{\bf Proof}$\ $    By Lemma \ref{E:harmonic}, Proposition \ref{II} and
the methods in Corollary 3.5 \cite{ROS}, we can prove (\ref{cd}) for
$k=0$.  The only difference is that we apply the estimate (\ref{O})
directly to replace  the scale invariant property of the Brownian
SLE used in \cite{ROS}.   Notice that for $k\geq 1$ and $t>0$, the
process $(V_{t-u}-V_{t-})_{0< u< t}$ has the same distribution with
$(V_u)_{0< u\leq t}$ conditional on $T_k=t$. From this fact and
$V_{t_2}=V_{t_2-} $ almost surely, we have
\begin{align}\label{flo}
&P\{|\hat{f}_{T_k}'(x-\triangle V_{T_k}+iy)|>\delta y^{-1}\}\nonumber\\
=& \int_{[t_1,t_2]}P\{|\hat{f}_{t}'(x-\triangle V_{T_k}+iy)|>\delta
y^{-1}| T_k=t\}P\{T_k\in dt\}
\nonumber\\
\leq & \int_{[t_1,t_2]}P\{|{f}_{t}'(x+V_{t-}+iy)|>\delta y^{-1} |
T_k=t
 \}P\{T_k\in dt\}
\nonumber\\
\leq & \sup_{{t_1}\leq t\leq {t_2}}P\{|{\hat{f}}_{t}'(x +iy)|>\delta
y^{-1} \}  .
\end{align}
By the result for  $k=0$,   (\ref{cd}) follows by (\ref{flo}). The
proof of (\ref{cde}) is similar to (\ref{cd}). See also Theorem 3.6
\cite{ROS}. \qed\medskip

\noindent By the  similar arguments as above, we can prove the
following result.
\begin{corollary}\label{dd21} Suppose that all the assumptions in
Proposition \ref{II} hold. Let $\sigma $ be a bounded nonnegative
random variable on $(\Omega,\mathcal{F},P)$ which is independent
with process $(V_t)_{t\geq 0}$. Then the estimate in (\ref{cde})
still holds with  $S_k$ replaced by $\sigma$.\end{corollary}

\begin{Lemma}\label{lemma}
Let $f$ be an analytic function on $\mathbb{H}$. Let $(r_j)_{j\geq
1}$ be a sequence of increasing positive numbers with $
\sum_{j=1}^\infty 1/r_j<\infty $. If for each $N\geq1$, there exists
a constant $C_4=C_4(N)$ such that
\begin{align}\label{98f7}
| {f}'(\frac{k}{2^j}+i2^{-j}) |&\leq C_42^{j}/r_j,\ \ \ \ \
k=0,\pm1,\cdots,\pm N2^j,
 \end{align}
  then $f$ can extend to
$\overline{\mathbb{H}}$ continuously.
\end{Lemma}\noindent{\bf Proof}$\ $ By distortion estimate (\ref{011}) and the condition (\ref{98f7}),
\begin{align}
\label{98f7f} | {f}'(x+iy) |&\leq C_4 K^32^{j}/r_j,\ \ \ \ \ |x|\leq
N, \ 2^{-j-1}<y\leq 2^{-j}.
\end{align}
Therefore for each $x\in [-N,N]$ and $y\in(0, 2^{-j}]$, we have
\begin{align}
|f(x+iy)-f(x+i2^{-j})|\leq C_4  K^3\sum_{n\geq j}1/r_n.
\end{align}
Hence for    $x_1,x_2\in [-N,N]$ with $|x_1-x_2|\leq 1/2^j$ and
$y_1,y_2\in(0, 2^{-j}]$,
\begin{align}
&|f(x_1+iy_1)-f(x_2+iy_2)|\nonumber\\\leq&
 |f(x_1+iy_1)-f(x_1+i2^{-j})|+|f(x_2+i2^{-j})-f(x_1+i2^{-j})|+
|f(x_2+iy_2)-f(x_2+i2^{-j})|\nonumber\\\leq& 3C_4\cdot
K^3\sum_{n\geq j}1/r_n,\nonumber
\end{align}
which leads to the conclusion by $ \sum_{j=1}^\infty 1/r_j<\infty
$.\qed \medskip

\begin{Lemma}\label{I61Idfff} Let $(V_t)_{t\geq 0}$ be the driving process in
(\ref{leq}) and let  $\sigma$ be a   nonnegative random variable on
$(\Omega,\mathcal{F},P)$. Suppose that  $\kappa\neq 4$ and  all the
parameters, such as $ \kappa_1,\kappa_2,\delta, b, a, \lambda,a'$,
satisfy the assumptions in Proposition \ref{II}. Suppose also that
the estimate $(\ref{cde})$ holds with $S_k$ replaced by
$\sigma\wedge M$ for any $M>0$ \emph{(}$ C_3$ in $(\ref{cde})$ may
depend on $M$\emph{)}. Then, almost surely, the conformal map
$f_{\sigma}(\cdot)$ extends to $\mathbb{\overline{H}}$ continuously.
\end{Lemma}\noindent{\bf Proof}$\ $ Case 1. $\kappa\in(0,4)\cup (4,8)$.
We assume that $\sigma\equiv t\geq 0$ because the proof  is similar.
Let $b=1\wedge [({4+\kappa})/{4\kappa}] $. By definition
(\ref{EQ:h_1:h_2}) we can check that for $\kappa_1-\kappa>0$ small
enough
 \begin{align}\label{odjf} \lambda-2b=2b+\kappa_1 b(1-2b)/2>1.\end{align}
  For any $\varepsilon>0$, we can find $M'>0$ such that
 \begin{align}\label{-}
 P\{A\}\geq1- \varepsilon,\ \ \ for\ A:=\{\sup_{0\leq u\leq t}|V_u|\leq
 M'\}.
 \end{align}
Let $N\geq 1$.  For each $j\geq1$, by estimate (\ref{cde}) with
$S_k$ replaced by $t$ and (\ref{odjf}), we can
 choose $\sigma>0$ small enough such that for some $\varepsilon'>0$
\begin{align}\label{anna11}
P\{|\hat{f}_{t}'(k2^{-j}+i2^{-j})|>  2^j2^{-\sigma
j}\}\leq& C_3(\kappa,\theta,b)(1+k^2)^b2^{-\lambda
j}2^{\lambda\sigma}\vartheta(2^{-\sigma j},a-\lambda) \nonumber\\
\leq &  C_3(\kappa,\theta,b)
 (1+(N+M')^2)^b 2^{-  j}2^{-\varepsilon' j},
  \end{align}
  for $k=0,\pm1,\cdots, \pm 2^j(N+M') $.
By (\ref{011}) and (\ref{anna11}), almost surely for $\omega\in A$
we can find $J(\omega)$ such that
$$
| {f}_{t}'(k/2^{-j}+i2^{-j})|\leq   e^2 2^j2^{-\sigma j},\ \ \
k=0,\pm1,\cdots, \pm 2^jN ,\ \ j\geq J(\omega).
$$
By Lemma \ref{lemma}, this implies that, almost surely for
$\omega\in A$, $f_{t}(\cdot)$   extends to $\mathbb{\overline{H}}$
continuously. Hence  the proof for this case is completed by taking
$\varepsilon\rightarrow0$.

Case 2. $\kappa\geq 8$. By taking  $b= ({4+\kappa})/{(4\kappa)} $ in
(\ref{EQ:h_1:h_2}), we can also check that   $
\lambda-2b=2b-\kappa_2 b(1-2b)/2>1 $ for $\kappa-\kappa_2>0$ small
enough. Therefore  we can prove the assertion as case 1.\qed\medskip

\begin{Remark} \label{rs} For $\kappa=4$,  we have  $\sup_{b>0}(\lambda-2b)=1$ and the maximum is equal to  $1$
at $b=1/2$, which is not enough to apply Lemma \ref{lemma}.
\end{Remark}

\noindent{\bf Proof of Theorem 1.1}$\ $  First we consider the
truncated case. Let $\delta_0$ be the positive constant  in
Proposition \ref{II} and let $(V_t) $ be the stochastic process
  defined by (\ref{gfffgg}) with $\delta\in (0,\delta_0)$. Choose $\beta$   such
that $1<\beta<2/\alpha$ and set  $L=[(2-\alpha\beta)^{-1}]+1$. For
$N\geq1$, $j\geq 1$ and $0\leq k\leq N2^{2j}-1$, define  stopping
times $(H_{j,k,l})_{1\leq l\leq L}$ by
 \begin{align}
H_{j,k,1}&= \inf \{ t\leq  ({k+1}){2^{-2j}}:
\theta^{1/\alpha}|\triangle S_{t } |>2^{-\beta j} \}\vee
{k}{2^{-2j}},\ \nonumber\\  \ H_{j,k,l}&=\inf \{ t<H_{j,k,l-1}:
\theta^{1/\alpha}|\triangle S_{t } |>2^{-\beta j} \}\vee
{k}{2^{-2j}},\ \ \ 2\leq l\leq L,\nonumber
 \end{align}    and define  stopping times $(R_{j,k,m})_{0\leq
m\leq j^2}$ by $R_{j,k,0}=(k+1)/2^{2j}$ and
\begin{align}
R_{j,k,1}&=  \inf \{ t<({k+1}){2^{-2j}}: |B_t-B_{{{(k+1)}}{2^{-2j}}}
|>2^{- j} \}\vee {k}{2^{-2j}},\ \ \nonumber\\ \  \ R_{j,k,m}&=\inf
\{ t<R_{j,k,m-1}: |B_t-B_{R_{j,k,m-1}} |>2^{- j} \}\vee
{k}{2^{-2j}},\ \ 2\leq m\leq j^2.\nonumber
 \end{align} By Lemma \ref{yes}, Lemma
 \ref{sure} and the distribution of the Brownian motion, we see
 that for a.s. $\omega\in \Omega$, there exists an integer
 $J_0=J_0(\omega)$ such that for $ j\geq J_0(\omega)$
\begin{align}\label{aa}
H_{j,k,L}(\omega)=R_{j,k, j^2}(\omega)&= {k}{2^{-2j}},\ \  \ 0\leq
k\leq N2^{2j}-1,
  \\
\sup_{k/2^{2j}\leq t< (k+1)/2^{2j}}|S_{2^{-\beta j},t}-S_{2^{-\beta
j},k/2^{2j}}|&\leq 2^{-j},\ \ \ \ \  \ 0\leq k\leq N2^{2j}-1,
  \\
\sup_{0\leq m\leq j^2-1}\sup_{R_{j,k,m+1}\leq t<
R_{j,k,m}}|B_t-B_{R_{j,k,m}}|&\leq 2^{-j},\ \ \ \ \ \ 0\leq k\leq
N2^{2j}-1.
\end{align}
 Let $\kappa\neq8$,
$b=((8+k)/4k)\wedge 1$ and assume that all the assumptions in
Proposition \ref{II} hold.
 By choosing $\kappa_1-\kappa$, $\kappa-\kappa_2$ small enough and applying
 Proposition \ref{oiu}, we can find $\varepsilon$ small
enough such that
\begin{align}  P\{|\hat{f}_{H_{j,k,l}}'(
i2^{-j})-\triangle V_{H_{j,k,l}}|>2^j/j^4\}\leq&
O(1)2^{-2j}2^{-\varepsilon j},\ \
 \ \nonumber \\\label{non}
 P\{\exists \ 1\leq m\leq j^2: |\hat{f}_{R_{j,k,
 m}}'(   i2^{-j})|>2^j /j^4\}\leq & O(1)2^{-2j}2^{-\varepsilon j}.\ \ \
\end{align}
 For   the proof of  $(\ref{non})$, see   (3.20) in \cite{ROS} for details. This
implies that almost surely,   there exists   an   integer
 $J_1=J_1(\omega)$ such that for  $j\geq J_1(\omega)$ and $0\leq k\leq N2^{2j}-1$
\begin{align}\label{anna11111}  |\hat{f}_{H_{j,k,l}}'(
i2^{-j})-\triangle V_{H_{j,k,l}}|&\leq 2^j /j^4 ,\ \ 1\leq l\leq
L;\\\label{aaaa}
  |\hat{f}_{R_{j,k,m}}'(   i2^{-j})|&\leq 2^j  /j^4,\ \  1\leq m\leq
 j^2.
\end{align}
Combing (\ref{aa})-(\ref{aaaa}) and applying Lemma \ref{lemma1} by
setting  $r_j=j^2$ and $c=2(\sqrt{\kappa}+\theta^{1/\alpha})$, we
can prove that $(\overline{\gamma}_t) $ is a c\'adl\'ag curve almost
surely, where $\overline{\gamma}_t$ is defined by (\ref{curve}) for
the  driving process $(V_t) $.

From now on  we further assume that $\kappa\neq4 $ and use  the
notations $\gamma_t, f_t,\zeta$,$K_t$   for the driving process
$(U_t) $ in (\ref{as}). As in Section 1, define $T_1$ to be the
first time $|\triangle S_t|>\delta$
  and  define by induction $T_{n+1}=\inf\{t\geq T_n: |\triangle S_t|>\delta\}$
  for $n\geq 1$. By  the arguments in Section 1, to obtain the  c\'adl\'ag property of
   $( {\gamma}_t) $ from that of $(\overline{\gamma}_t) $, we only need to prove that
   $f_{T_k},k\geq1$, can extend to $\mathbb{\overline{H}}$ continuously. For $k=1$, noticing    that
   $U_t=V_t$ for $t\in [0,T_1)$ and $T_1$ is
   independent with $(V_t) $, this follows by
   Proposition \ref{oiu}, Corollary  \ref{dd21} and Lemma \ref{I61Idfff}.
Suppose by induction that it holds for $k=n$. For $k=n+1$, by
  relation $f_{T_{n+1}}=f_{T_n}(g_{T_n}\circ f_{T_{n+1}})$, the
assertion for $f_{T_{n+1}}$ can reduce  to the continuous extension
of conformal map $g_{T_n}\circ f_{T_{n+1}}$. By considering the
truncated process
$$U_{T_n}+\sqrt{\kappa}B_{t+T_n}+\theta^{1/\alpha}S_{\delta,t+T_n}
-\sqrt{\kappa}B_{ T_n}-\theta^{1/\alpha}S_{\delta, T_n},\ \ \ t\geq
0,$$ the continuous extension of   $g_{T_n}\circ f_{T_{n+1}}$ can be
proved as the case $k=1$. Hence we obtain  the first assertion of
the theorem.

Next  we turn to the relation between $(K_t)_{t\geq0}$ and
$(\gamma_t)_{t\geq 0}$. Denote $$A=\{w: (\gamma_t)_{t\geq 0} \mbox{\
is a c\'adl\'ag curve}\}.$$  We have $P\{A\}=1$. Let $B_t$ be the
unbounded component of $\mathbb{H}\setminus \gamma[0,t]$. We want to
show that $B_t=H_t:=\mathbb{H}\setminus K_t$ for $\omega\in A$. If
this is not true, then   there exist some $\omega\in A$ and $t> 0$
such that $H_t \neq B_t$. By that $H_t\subseteq B_t$ and that
$\gamma[0,t] $ is closed, we can find $z_0\in
\partial K_t\setminus \gamma[0,t]$ with
 $dist(z_0,\gamma[0,t]) >0$. Choose   $z_1\in H_{\zeta(z_0)}$ with $|z_1-z_0|<dist(z_0,\gamma[0,t])$ and
set $z_2=z_1+s_0(z_0-z_1)$ with $s_0=\inf\{s>0:z_1+s(z_0-z_1)\in
\partial K_{\zeta({z_0})}\}$. By  the definition of
$\zeta(z_2)$ and Proposition 2.14 \cite{CH}, we have     $\lim_{t
\uparrow 1}g_{\zeta({z_2})}(\xi_t)=U_{\zeta(z_2)}$, where
$\xi(t)=(1-t)z_1+tz_2$. Set
$\overline{\xi}_t=g_{\zeta({z_2})}^{-1}(i(1-t))$ for $0<t<1$, we
also have $\lim_{t \uparrow
1}g_{\zeta({z_2})}(\overline{\xi}_t)=U_{\zeta(z_2)}$. Therefore we
have $z_2=\gamma(\zeta(z_2))$ by Proposition 2.14 \cite{CH}, which
contradicts with $dist(z_2,\gamma[0,t])>0$.

To complete the proof of Theorem 1.1 for $\kappa \neq 4,8$, we only
need to prove that   almost surely, $f_{t}(\cdot)$ extends to
$\mathbb{\overline{H}}$ continuously for all  $t\in[0,\infty)$.   By
Theorem 2.1  \cite{CH}, for any $t\geq 0$, the continuous extension
of $f_t$ is equivalent to the local connected property of $\partial
K_t$. Thus by Proposition \ref{oiu}  and Lemma \ref{I61Idfff},
$\partial K_t$ is local connected for $t\in Q_+$ almost surely,
where $Q_+$ is the set of positive rational numbers. Since $\partial
K_t \subseteq \gamma([0,t])$ by $B_t=H_t$ and $(\gamma_t) $ is right
continuous, we have the following right continuity of $\partial K_t$
\begin{align}\label{conti}
\lim_{h\downarrow0}\sup\{|x-y|:x,y\in \partial K_{t+h}\setminus
\partial K_t \}=0, \ \ \ t\geq 0.
\end{align}
This can extend the locally connected property of   $\partial K_t$ 
from  $t\in Q_+$ to  $t\in [0,\infty)$.  Thus   by the equivalent conditions  in Theorem 2.1
\cite{CH} once  again, the continuous extension of $f_t$ holds for
all $t\in[0,\infty)$ almost surely.
   \qed\medskip

 At last we give an example to show that some  non locally connected sets    can be generated by a c\`ad\`ag curve.   Define a comb space by  $$D=\cup_{n=1}^\infty \{(x,y): x=1/n,\ 0\leq y\leq
   1\}\cup \{(x,y):x=0,\ 0\leq y\leq 1\}.$$
   Next we show that there exists a c\`ad\`ag curve  $(\xi(t))_{0\leq t\leq 2}$ in $\R^2$
   such that \begin{align}\label{dy}
D=\xi[0,2]:=\overline{\{\xi(t):0\leq t\leq 2\}},
   \end{align} and $\xi[0,t]$ is a connected set in $\R^2$ for each
   $t\in[0,2]$. For two subsets $A$ and $B$ of $\R$, we write $A+B=\{x+y:x\in A,y\in
   B\}$.
Set \begin{align}
R_1=&\{1\}+[1/4,1/2)\cup [3/4,1):=I_{1,1}\cup I_{1,2},\nonumber\\
R_2=&\{1\}+[1/16,2/16)\cup[3/16,4/16)\cup[9/16,10/16)\cup[11/16,12/16):=I_{2,1}\cup
I_{2,2}\cup I_{2,3}\cup I_{2,4}.\nonumber\end{align} Inductively,
define for $n\geq 2$
\begin{align}R_{n+1}= \{1\}+\big([0,1) \backslash R_{n} \big)\backslash
\bigcup_{k=0}^{2^{2n+1}}[(2k)2^{-2(n+1)},(2k+1)2^{-2(n+1)}):
=\bigcup_{k=1}^{2^{n+1}}I_{n+1,k},
\end{align}
where $(I_{n+1,k})$ are intervals of length $2^{-2(n+1)}$ arranged
by the  increasing order. Set  $\xi(t)=(t,0)$ for $0\leq t<1$ and
$\xi(2)=(0,0)$. Define for $n\geq 1$
 and $1\leq k\leq 2^{n}$\begin{align}
\xi_t= (1/n, (k-1)2^{-n}+(t-t_{n,k})2^{n}),\ \ \ t\in I_{n,k},
\end{align}
where $t_{n,k}=\inf I_{n,k}$. For each point $t\in[1,2)\backslash
\cup_{n\geq 1}R_n$, we can find a sequence of integers $(k_n)$ such
that $(t_{n,k_n})$ is a decreasing sequence in
$\cup_{n\geq 1}R_n$ and
$\lim_{n\rightarrow\infty}t_{n,k_n}=t$. By   definition, it is
easy to check that $\xi({t_{n,k_n}})$ converges to a number which does not
depend on the choice of $(t_{n,k_n})$. Define
\begin{align}
\xi(t)=(0,\lim_{n\rightarrow \infty}\xi({t_{n,k_n}})).
\end{align}
By definition we can check directly that $\xi(t)$ is a right
continuous and left limit curve satisfying (\ref{dy}). Notice  that
for each $n= 1,2\cdots, \infty$ and $0\leq y_1<y_2\leq 1$,  we have $\inf \{s: (1/n,y_1\in
\xi[0,s]\}<\inf \{s:(1/n, y_2)\in \xi[0,s]\}$, where $1/\infty =0$. This implies that
$\xi[0,t]$ is   connected   for each
   $t\in[0,2]$.

 \subsection*{Acknowledgements} The author thanks the supports of  EPSRC grant
GR/T26368/01 for the starting of this work in the Statistics
Department
 of Oxford University and also  the  helpful  discussions with Matthias Winkel on this topic.
 The author   thanks  Zhen-Qing Chen for the  helpful comments.

 \end{document}